\documentclass[11pt]{amsart}
\usepackage{amsmath,amssymb} 
\usepackage{graphicx}
\usepackage[utf8]{inputenc} 
\usepackage[font=small,labelfont=bf]{caption}

\begin{document}

\newtheorem{thm}{Theorem}[section]
\newtheorem{lem}[thm]{Lemma}
\newtheorem{prop}[thm]{Proposition}
\newtheorem{coro}[thm]{Corollary}
\newtheorem{defn}[thm]{Definition}
\newtheorem{remark}{Remark}

\numberwithin{equation}{section}

\newcommand{\Z}{{\mathbb Z}} 
\newcommand{\Q}{{\mathbb Q}}
\newcommand{\PP}{{\mathbb P}}
\newcommand{\R}{{\mathbb R}}
\newcommand{\C}{{\mathbb C}}
\newcommand{\N}{{\mathbb N}}
\newcommand{\FF}{{\mathbb F}}
\newcommand{\T}{{\mathbb T}}
\newcommand{\fq}{\mathbb{F}_q}

\newcommand{\fixmehidden}[1]{}

\def\scrA{{\mathcal A}}
\def\cB{{\mathcal B}}
\def\Eps{{\mathcal E}}
\def\cI{{\mathcal I}}
\def\scrD{{\mathcal D}}
\def\cF{{\mathcal F}}
\def\cL{{\mathcal L}}
\def\cM{{\mathcal M}}
\def\cN{{\mathcal N}}
\def\cP{{\mathcal P}}
\def\scrR{{\mathcal R}}
\def\scrS{{\mathcal S}}

\newcommand{\rmk}[1]{\footnote{{\bf Comment:} #1}}

\renewcommand{\mod}{\;\operatorname{mod}}
\newcommand{\ord}{\operatorname{ord}}
\newcommand{\TT}{\mathbb{T}}
\renewcommand{\i}{{\mathrm{i}}}
\renewcommand{\d}{{\mathrm{d}}}
\renewcommand{\^}{\widehat}
\newcommand{\HH}{\mathbb H}
\newcommand{\Vol}{\operatorname{vol}}
\newcommand{\area}{\operatorname{area}}
\newcommand{\tr}{\operatorname{tr}}
\newcommand{\norm}{\mathcal N} 
\newcommand{\intinf}{\int_{-\infty}^\infty}
\newcommand{\ave}[1]{\left\langle#1\right\rangle} 
\newcommand{\E}{\mathbb E}
\newcommand{\Var}{\operatorname{Var}}
\newcommand{\Cov}{\operatorname{Cov}}
\newcommand{\Prob}{\operatorname{Prob}}
\newcommand{\sym}{\operatorname{Sym}}
\newcommand{\disc}{\operatorname{disc}}
\newcommand{\CA}{{\mathcal C}_A}
\newcommand{\cond}{\operatorname{cond}} 
\newcommand{\lcm}{\operatorname{lcm}}
\newcommand{\Kl}{\operatorname{Kl}} 
\newcommand{\leg}[2]{\left( \frac{#1}{#2} \right)}  
\newcommand{\id}{\operatorname{id}}
\newcommand{\beq}{\begin{equation}}
\newcommand{\eeq}{\end{equation}}
\newcommand{\bsp}{\begin{split}}
\newcommand{\esp}{\end{split}}
\newcommand{\bra}{\left\langle}
\newcommand{\ket}{\right\rangle}
\newcommand{\diam}{\operatorname{diam}}
\newcommand{\supp}{\operatorname{supp}}
\newcommand{\dist}{\operatorname{dist}}
\newcommand{\sgn}{\operatorname{sgn}}
\newcommand{\inte}{\operatorname{int}}
\newcommand{\Spec}{\operatorname{Spec}}
\newcommand{\ddiv}{\operatorname{div}}
\newcommand{\sumstar}{\sideset \and^{*} \to \sum}

\newcommand{\LL}{\mathcal L} 
\newcommand{\sumf}{\sum^\flat}
\newcommand{\Hgev}{\mathcal H_{2g+2,q}}
\newcommand{\USp}{\operatorname{USp}}
\newcommand{\conv}{*}
\newcommand{\CF}{c_0} 
\newcommand{\kerp}{\mathcal K}

\newcommand{\gp}{\operatorname{gp}}
\newcommand{\Area}{\operatorname{Area}} 


\title[Multifractal scaling and the Euler equations]
{Multifractal scaling and \\ the Euler equations on $\R^3/\Z^3$}

\author{Henrik Uebersch\"ar}
\address{Sorbonne Universit\'e, Universit\'e Paris Cit\'e, CNRS, IMJ-PRG, F-75006 Paris, France.}
\email{henrik.ueberschar@imj-prg.fr} 
\date{\today} 

\maketitle 
 
\begin{abstract}
We study the Euler equations describing the motion of an incompressible fluid on the cubic torus with real initial data. We construct solutions on the Fourier side which display a sudden loss of regularity within finite time even for highly regular initial data. Moreover, the solution may regain its initial regularity within finite time. This loss of regularity may coincide with the appearance of a certain type of multifractal scaling of the solutions. 
\end{abstract}

\section{Introduction}
In the present article, we investigate the occurrence of multifractal scaling for certain solutions to the incompressible Euler equations on the cubic torus $\T^3=\R^3/\Z^3$. We present an elementary construction of solutions on the Fourier side which may evolve from arbitrarily high initial regularity to arbitrarily bad behaviour along at least one direction of the lattice $\Z^3$ within finite time. We apply this result to construct solutions, where this loss of regularity coincides with a certain type of multifractal scaling, which is explained in section \ref{sec-mult}.

One particular feature of our solutions is their non-uniqueness on the time-interval on which the loss of regularity and the multifractal scaling occurs. While a well-behaved solution exists, there exists a multitude of other solutions which display irregular behaviour on the respective time-interval. The non-uniqueness of the Euler equations has been known since at least 1993, when Scheffer proved the existence of weak solutions on the plane which were compactly supported in space and time \cite{Sch93}. 
In 1996, Snirelman \cite{Sn96} gave a simpler construction for compact domains. While both Scheffer's and Snirelman's constructions were quite complicated, in 2009, de Lellis and Szekelyhidi gave an alternative construction by elementary methods \cite{dLSz09}. In 2013, the same authors constructed weak solutions on $\T^3$ which dissipate the kinetic energy completely \cite{dLSz13}.


\subsection{Set-up.}
Let $u:\T^3\times\R_+ \to \R^3$ and $p:\T^3\times\R_+\to \R$ denote a solution to the Euler system, where external forcing is absent:
\beq\label{Euler}
\begin{split}
\begin{cases}
\partial_t u&=-(u\cdot\nabla)u-\nabla p,\\
\\
\ddiv u&=0.\\
\\ 
u(\cdot,0)&=u_0
\end{cases}
\end{split}
\eeq
where $u_0$ denotes the initial data, a divergence-free vector field on $\T^3$.


We define the Fourier mode for each component of the vector field $u=(u_1,u_2,u_3)$:
$$\hat{u}_j(\xi,t):=\int_{\T^3}u_j(x,t) e_{-\xi}(x) dx, \quad\text{where}\quad e_\xi(x):=e^{2\pi i \xi\cdot x}.$$

We will study the Euler system on the Fourier side by looking for a solution $\hat{u}=(\hat{u}_1,\hat{u}_2,\hat{u}_3):\Z^3\mapsto\C^3$ to the system of coupled differential equations for any $\xi\in\Z^3$ which arises by applying the Fourier transform to \eqref{Euler}. We note that a real solution must satisfy the symmetry relation $\hat{u}(-\xi,\cdot)=\overline{\hat{u}(\xi,\cdot)}$. Moreover, we take $p(x,t)=p_0(t)$ to be constant with respect to the space variable and we assume $p_0\in C^\infty(\R_+)$. 

This gives rise to the system: 
\beq\label{Fourier-system}
\begin{cases}
\begin{split}
\partial_t \hat{u}(\xi,\cdot)
=&-2\pi i\sum_{\substack{\zeta,\eta\in\Z^3 \\ \zeta+\eta=\xi}}(\hat{u}(\zeta,\cdot)\cdot\eta)\hat{u}(\eta,\cdot),\\
\\
\hat{u}(\xi,\cdot)\cdot\xi=&0,\\
\\
\hat{u}(\xi,0)=&\hat{u}_0(\xi).
\end{split}
\end{cases}
\eeq
The advantage of studying the system \eqref{Fourier-system} is that this enables us to track the regularity of any solution in terms of the decay properties of the norm $\|\hat{u}(\xi,t)\|$, as $|\xi|\to\infty$, as time evolves. It is known \cite{FFT75}, that in dimension $d=2$, for the case of real-valued solutions, a certain decay of the Fourier modes of the initial data on the lattice may be preserved for all $t>0$. This is also known for the related case of the Navier-Stokes system \cite{FT89},\cite{SM99}. However, the situation may be different if one admits complex solutions. For instance, Sinai and Li proved blow-up for complex solutions of the Navier-Stokes system \cite{CS17}.

\subsection{Main result}
We consider certain real initial data and show that there exists a multitude of solutions. While there is a solution which preserves regularity, there are also many other solutions which, within finite time, evolve to an arbitrarily bad behaviour on $\Z^3$ (sections \ref{sec-der}, \ref{sec-soln}) in the following sense. 
\begin{thm}\label{3d-blowup}
Fix two functions $f_1,f_2:\Z^3\to\R_+$.
There exists a solution to \eqref{Fourier-system} such that the initial data satisfy $|u_0(\xi)|\leq f_1(\xi)$ for $|\xi|>1$, the $\hat{u}(\xi,\cdot)$ are smooth in the time variable 
and there exists a time $T_0>0$ and a vector $\xi_1\in\Z^3\setminus\{0\}$ such that $|\hat{u}(m\xi_1,T_0)|=f_2(m\xi_1)$ for any $m\in\Z\setminus\{0\}$.

\end{thm}\label{blow-up}
\begin{remark}
This theorem can, for example, be applied to initial data satisfying an exponential decay $f_1(\xi)=e^{-|\xi|}$ and for an algebraic function $f_2(\xi)=|\xi|^{\alpha}$, at time $T_0>0$, along the direction $\xi=m\xi_1$, $m\in\Z\setminus\{0\}$, for a real exponent $\alpha$. In particular, we may blow up any Sobolev norm within finite time (section \ref{sec-blowup}), and we may use this to show that the solution (or rather certain transversal sections of it) display multifractal scaling, as soon as the solution loses regularity (section \ref{sec-mult}).
\end{remark}
\begin{remark}
In order to construct these initial data and solutions it is essential to work in dimension $3$. However, if one admits {\em complex-valued} initial data, it is possible to construct a solution which loses regularity even in dimension $2$. This is discussed in section \ref{sec-complex}.
\end{remark}
 
We point out that it is easily possible to modify our construction to give solutions which suddenly lose regularity on a specified time-interval to regain the original regularity at a later time, and, in fact, this can be repeated infinitely many times (section \ref{sec-gain}). Finally, we note that the type of construction presented here may be applied to more general non-linear PDE which is the subject of a forthcoming article \cite{U}.

\section{Multifractal scaling}\label{sec-mult}
A simple consequence of Theorem \ref{blow-up} is that the loss of regularity of the solution $u$, as time evolves, may coincide with a certain type of multifractal scaling of the solution. Recently, the author proved, in joint work with Keating, the existence of multifractal eigenfunctions for certain models in quantum chaos, whose dynamics is intermediate between chaos and integrability: singular quantum billiards \cite{KU1} and quantum star graphs \cite{KU2}. In both cases, the method of proof relied on the explicit nature of these models, as well as, in the first case, certain arithmetic symmetries satisfied by the  quantum billiard.

In this section, we will define the Renyi entropy associated with certain measures on the circle 
$\T^1=\R/\Z$ which arise by truncating the Fourier series and integrating out transversal sections of each component of the vector field. This gives rise to a vector field on the circle. If we apply the Euclidean norm at each point, then we obtain a square integrable function on the circle. Upon $L^2$-normalizing this function, we obtain a probability measure on $\Z$ by Plancherel's theorem.

Denote $x:=(x_1,x_2,x_3)\in\T^3$. We have the Fourier series
$$
u(x_1,x_2,x_3;t)=\sum_{\xi\in\Z^3}\hat{u}(\xi,t)e_\xi(x_1,x_2,x_3).
$$
Consider an approximation to an actual solution $u$ which is given by a cutoff in Fourier space (take $N>0$)
$$
u_N(x_1,x_2,x_3;t)=\sum_{|\xi|\leq N}\hat{u}(\xi,t)e_\xi(x_1,x_2,x_3).
$$

Let us consider the vector field $U_N:\T^1\to\R^3$ given by integrating out two space variables (recall our choice $\xi_1=(1,0,0)$):
$$
(U_N)_j(x_1,t)=\int_{\T^2}(u_N)_j(x_1,x_2,x_3;t)dy\,dz=\sum_{m\in\Z} (\hat{u}_N)_j(m\xi_1,t)e^{2\pi i mx_1}
$$
Hence, $(\hat{U}_N)_j(m,t)=\hat{u}_j(m\xi_1,t)$.

We will study the following discrete probability measure on $\Z$: $$\mu_N(m,t)=\frac{|\hat{U}_N(m\xi_1,t)|^2}{\sum_{m'\in\Z}|\hat{U}_N(m'\xi_1,t)|^2}.$$

We will see that Theorem \ref{blow-up} allows us to construct solutions which lead to a multifractal scaling law, which is a power law for the qth moments of the measure $\mu_N$ of the form $$\sum_{|m|\leq N}\mu_N(m,t)^q \sim N^{(1-q)D_q(t)},$$ where the scaling exponent $D_q$ varies as a function of $q$. 

This motivates our definition of the Renyi entropy associated with the measure $\mu_N(\cdot,t)$: 
\begin{equation}
H_{q,N}(t):=\frac{1}{1-q}\log\left(\sum_{|m|\leq N}\mu_N(m,t)^q\right)
\end{equation}
If we can obtain an asymptotic of the form $H_{q,N}(t)\sim c_q \log N$, as $N\to+\infty$, then we may compute the scaling exponent explicitly as
\begin{equation}
D_q(t)=\lim_{N\to+\infty} \frac{H_{q,N}(t)}{\log N}.
\end{equation}

We have the following corollary of Theorem \ref{blow-up} which we prove in section \ref{sec-mult-proof}.
\begin{coro}\label{mult}
There exist initial data and a time $T_0$ as in Theorem \ref{blow-up} with $f_2(\xi)=|\xi|^{-\alpha}$, $\alpha\in(0,1/2)$, such that for any $t>T_0$ and any $q>\tfrac{1}{2\alpha}$ we have
\begin{equation}
D_q(t)=(1-2\alpha)\frac{q}{q-1}.
\end{equation}
\end{coro}
 
\section{Real data on the $3$-torus}

\subsection{Derivation}\label{sec-der}
Let us consider solutions $u:\T^3\times\R_+\mapsto\R^3$. We fix a vector $v\in\Z^3$ and denote $S=\bra v\ket^\perp\cap\Z^3$.
\footnote{We use the notation $\bra v\ket=\{\lambda v \mid \lambda\in\R\}$.} Choose $\xi_0,\xi_1,\eta_0\in S\setminus\{0\}$ such that $\eta_0\cdot\xi_1=0$. Throughout this section we will take, for simplicity, $v=(0,0,1)$, $\eta_0=\xi_0=(0,1,0)$ and $\xi_1=(1,0,0)$. Hence, $S=\Z\eta_0\oplus\Z\xi_1$.

Let us look for a solution such that $\hat{u}:\Z^3\times\R_+\to\C$ vanishes outside $S$ (and we must impose the symmetry $\hat{u}(\xi,\cdot)=\overline{\hat{u}(-\xi,\cdot)}$, since we are looking for real-valued solutions). 

We assume that $\hat{u}(\xi,\cdot)\in\bra v\ket$ except for $\xi=0,\eta_0,-\eta_0$, where we set $\hat{u}(0,\cdot)\in\bra \xi_0\ket$ and $\hat{u}(\pm \eta_0,\cdot)\in\bra \xi_1\ket$. In particular, $\hat{u}(0,0)=\xi_0$ and $\hat{u}(\pm\eta_0,0)=\xi_1$.
We verify that the vector field is divergence-free (recall $\eta_0\cdot\xi_1=0$).

We recall the differential equation for the Fourier modes in \eqref{Fourier-system}
\begin{equation}\label{Fourier-DE}
\begin{cases}
\begin{split}
\partial_t \hat{u}(\xi,\cdot)
=&-2\pi i\sum_{\substack{\zeta,\eta\in\Z^3 \\ \zeta+\eta=\xi}}(\hat{u}(\zeta,\cdot)\cdot\eta)\hat{u}(\eta,\cdot),\\
\\
\hat{u}(\xi,\cdot)\cdot\xi=&0,\\
\\
\hat{u}(\xi,0)=&\hat{u}_0(\xi).
\end{split}
\end{cases}
\end{equation}

{\bf Case1: $\xi\notin S$.} Because $\eta+\zeta=\xi$, it follows that either $\eta\notin S$ or $\zeta\notin S$. Because of our assumption that $\hat u$ vanishes outside $S$ it follows that all terms in the sum must vanish. Hence $\partial_t \hat{u}(\xi,\cdot)=0$ and $\hat{u}(\xi,t)=\hat{u}(\xi,0)$.

{\bf Case 2: $\xi\in S$.} In this case we have, since $\hat{u}(\eta,\cdot)$ vanishes outside $S$, that the sum runs over terms such that
 $\zeta=\xi-\eta\in S$. Note that $\hat{u}(\zeta,\cdot)\cdot\eta=0$ unless $\zeta=0,\eta_0,-\eta_0$, because, in this case, $\hat{u}(\zeta,\cdot)\in\bra v\ket$ and $\eta\in S$. So only $3$ terms on the r.h.s. of \eqref{Fourier-DE} contribute:
\begin{equation}\label{Fourier-DE-2}
\begin{split}
\partial_t \hat{u}(\xi,\cdot)
=&-2\pi i((\hat{u}(0,\cdot)\cdot\xi)\hat{u}(\xi,\cdot)\\
&\quad+(\hat{u}(\eta_0,\cdot)\cdot(\xi-\eta_0))\hat{u}(\xi-\eta_0,\cdot)\\
&\quad+(\hat{u}(-\eta_0,\cdot)\cdot(\xi+\eta_0))\hat{u}(\xi+\eta_0,\cdot))\\
\end{split}
\end{equation}
which we simplify as
\begin{equation}
\begin{split}
\partial_t \hat{u}(\xi,\cdot)
=&-2\pi i((\hat{u}(0,\cdot)\cdot\xi)\hat{u}(\xi,\cdot)\\
&\quad+(\hat{u}(\eta_0,\cdot)\cdot\xi)\hat{u}(\xi-\eta_0,\cdot)\\
&\quad+(\hat{u}(-\eta_0,\cdot)\cdot\xi)\hat{u}(\xi+\eta_0,\cdot))
\end{split}
\end{equation}

In particular, we have 
$$\partial_t \hat{u}(k\eta_0,\cdot)=-2\pi i((\hat{u}(0,\cdot)\cdot k\eta_0)\hat{u}(k\eta_0,\cdot)$$
because $\hat{u}(\eta_0)\cdot\eta_0=0$. So if $k=0$: $\partial_t \hat{u}_j(0,\cdot)=0$ and therefore $\hat{u}(0,t)=\hat{u}(0,0)=\xi_0$. We obtain  
\begin{equation}
\hat{u}(k\eta_0,t)=e^{-2\pi i k(\xi_0\cdot \eta_0)t}\hat{u}(k\eta_0,0).
\end{equation}

So, we may rewrite \eqref{Fourier-DE-2} as
\begin{equation}
\partial_t \hat{u}(\xi,\cdot)=-2\pi i[(\xi_0\cdot\xi)\hat{u}(\xi,\cdot)+e^{-2\pi i(\xi_0\cdot \eta_0)t}(\xi_1\cdot\xi)\hat{u}(\xi-\eta_0,\cdot)+e^{2\pi i(\xi_0\cdot \eta_0)t}(\xi_1\cdot\xi)\hat{u}(\xi+\eta_0,\cdot)]
\end{equation}

Recall our choice $v=(0,0,1)$, $\eta_0=(0,1,0)=\xi_0$ and $\xi_1=(1,0,0)$, as well as $S=\Z\eta_0\oplus\Z\xi_1$.
We set $v_{k,m}(t)=\hat{u}(k\eta_0+m\xi_1,t)$, $k\in\Z$, $m\in\Z\setminus\{0\}$, and rewrite the differential equation as
$$v_{k,m}'=\alpha_m(t)v_{k-1,m}+\beta_k v_{k,m}+\widetilde\alpha_m(t)v_{k+1,m},$$ where $\alpha_m(t)=-2\pi i m e^{-2\pi it}$, $\widetilde\alpha_m(t)=-\overline{\alpha_m(t)}$ and $\beta_k=-2\pi i k$.

Now, we rewrite to obtain the recurrence relations
\begin{equation}
v_{k+1,m}=-\frac{1}{\overline{\alpha_m(t)}}v_{k,m}'+\frac{\alpha_m(t)}{\overline{\alpha_m(t)}}v_{k-1,m}
+\frac{\beta_k}{\overline{\alpha_m(t)}}v_{k,m}
\end{equation}
and
\begin{equation}
v_{k-1,m}=\frac{1}{\alpha_m(t)}v_{k,m}'+\frac{\overline{\alpha_m(t)}}{\alpha_m(t)}v_{k+1,m}-\frac{\beta_k}{\alpha_m(t)} v_{k,m}.
\end{equation}

So, for given functions $v_{0,m}$ and $v_{1,m}$ we may determine $v_{k,m}$ for any $k\in\Z$ using the above recurrence relations.

Fix $T>0$. Let us introduce a function $f\in C^\infty(\R_+)$ which is supported away from $[0,T]$. For example, we may take
\begin{equation}
f(t)=
\begin{cases}
0, \quad \text{if}\;t\leq T,\\
\\
\exp(\frac{1}{T-t}), \quad \text{if}\; t>T.
\end{cases}
\end{equation}

Moreover, set $v_{0,m}=g(m)f v$, for a given $g:\Z\setminus\{0\}\to\C$, which satisfies $g(-m)=\overline{g(m)}$, and set $v_{1,m}(t)=0$ for all $t\geq0$. As a result of the recurrence relation we see that for any $k\in\Z$ we have $\supp v_{k,m}\subset[T,+\infty)$, because for any $n\geq0$: $\supp f^{(n)}\subset[T,+\infty)$, and $v_{k,m}\in\bra v\ket$.

We have, thus, derived a solution which satisfies
\beq
\begin{split}
\hat{u}(m\xi_1+k\eta_0,t)=
\begin{cases}
e^{-2\pi i(\xi_0\cdot\eta_0)t}\hat{u}(k\eta_0,0), \quad \text{if}\; m=0,\\
\\
v_{k,m}(t), \quad\text{if}\; m\neq 0,
\end{cases}
\end{split}
\eeq
and $\hat{u}(0,0)=\xi_0$, $\hat{u}(\pm\eta_0,0)=\xi_1$, whereas for $|k|>1$: $\hat{u}(k\eta_0)=h(k)v$, and $h:\Z\setminus\{-1,0,1\}\to \C$ satisfies $h(-k)=\overline{h(k)}$.

\subsection{Construction of the solution}\label{sec-soln}
First of all, we set $\hat{u}(\xi,t)=0$ for all $\xi\notin S$ and for all $t\geq0$.

We define 
\begin{equation}
\hat{u}(k\eta_0,t)=
\begin{cases}
\begin{split}
\xi_0, &\qquad \text{if\;} k=0,\\
e^{-2\pi i k t}\xi_1, &\qquad \text{if\;} k=1,-1,\\
e^{-2\pi i k t}h(k)v, &\qquad \text{if\;} |k|>1.
\end{split}
\end{cases}
\end{equation}

We set $\hat{u}(k\eta_0+m\xi_1,t)=v_{k,m}(t)$ for $m\in\Z\setminus\{0\}$ and $k\in\Z$.

Let
\begin{equation}
v_{0,m}=g(m)f(t)v, \qquad m\in\Z\setminus\{0\},
\end{equation}
and
\begin{equation}
v_{1,m}=v_{-1,-m}=0, \qquad m\geq 1.
\end{equation}

To define the Fourier modes for the remaining lattice vectors in $S$ we set $\alpha_m(t)=-2\pi i m e^{-2\pi it}$ and $\beta_k=-2\pi i k$.

Now, we define the following recurrence relations
\begin{equation}\label{R1}
v_{k+1,m}=-\frac{1}{\overline{\alpha_m}}v_{k,m}'+\frac{\alpha_m}{\overline{\alpha_m}}v_{k-1,m}
+\frac{\beta_k}{\overline{\alpha_m}}v_{k,m}
\end{equation}
and
\begin{equation}\label{R2}
v_{k-1,m}=\frac{1}{\alpha_m}v_{k,m}'+\frac{\overline{\alpha_m}}{\alpha_m}v_{k+1,m}-\frac{\beta_k}{\alpha_m} v_{k,m}.
\end{equation}

So for any given $m\in\Z\setminus\{0\}$, we have for $m\geq 1$ the pair of Fourier modes $(v_{0,m},v_{1,m})$ which then determines $v_{k,m}$ for all $k\in\Z$, and for $m\leq-1$ the pair $(v_{0,m},v_{-1,m})$ determines $v_{k,m}$ for all $k\in\Z$. Because all derivatives of $f$ are supported away from $[0,T]$ we have that $v_{k,m}$ is supported away from $[0,T]$.

We claim that $\hat{u}:\Z^3\times\R_+\to\C$ is a solution of \eqref{Fourier-system} with initial condition 
\begin{equation}
\hat{u}_0(\xi)=
\begin{cases}
\begin{split}
\xi_0, &\quad\text{if}\;\xi=0,\\
\xi_1, &\quad\text{if}\;\xi=\pm\eta_0,\\
h(k)v, &\quad\text{if}\;\xi=k\eta_0, \;|k|>1,\\
0, &\quad\text{otherwise.}
\end{split}
\end{cases}
\end{equation}
In particular, since $h(-k)=\overline{h(k)}$, we see that $u_0$ is real. 

Moreover, the recurrence relations \eqref{R1}, \eqref{R2} satisfy symmetry relations which ensure that the solution is real: if we take the complex conjugate of \eqref{R1}, then we obtain (using $v_{-k,-m}=\overline{v_{k,m}}$ as well as $\alpha_{-m}=-\alpha_m$ and $\beta_{-k}=\overline{\beta_k}$)
$$v_{-k-1,-m}=\frac{1}{\alpha_{-m}}v_{-k,-m}'+\frac{\overline{\alpha_{-m}}}{\alpha_{-m}}v_{-k+1,-m}-\frac{\beta_{-k}}{\alpha_{-m}} v_{-k,-m}.$$

Because $v_{0,m}\in\bra v\ket$, the recurrence relation ensures $v_{k,m}\in \bra v\ket$ for all $k\in\Z$ and $m\in\Z\setminus\{0\}$. Hence, we see that for all $\xi\in S\setminus\{0,\pm\eta_0\}$ we have $\hat{u}(\xi,t)\in \bra v\ket$. So $\hat{u}(\xi,t)\cdot\xi=0$ for all $\xi\in\Z^3$ and all $t\geq0$, because $\hat{u}(\pm\eta_0,t)\in\bra\xi_1\ket$ and $\xi_1\cdot\eta_0=0$.

To see that the differential equation 
$$\partial_t \hat{u}(\xi,\cdot)
=-2\pi i\sum_{\substack{\zeta,\eta\in\Z^3 \\ \zeta+\eta=\xi}}(\hat{u}(\zeta,\cdot)\cdot\eta)\hat{u}(\eta,\cdot)$$
is satisfied we need to verify this for any $\xi\in S$.

First of all, because $\hat{u}$ vanishes away from $S$ we see that the sum on the r.h.s. only runs over lattice vectors in $S$. The sum then simplifies 
\begin{equation}
\begin{split}
\sum_{\substack{\zeta,\eta\in S \\ \zeta+\eta=\xi}}(\hat{u}(\zeta,\cdot)\cdot\eta)\hat{u}(\eta,\cdot)
=&\sum_{\zeta=0,\pm\eta_0}(\hat{u}(\zeta,\cdot)\cdot(\xi-\zeta))\hat{u}(\xi-\zeta,\cdot)\\
=&\sum_{\zeta=0,\pm\eta_0}(\hat{u}(\zeta,\cdot)\cdot\xi)\hat{u}(\xi-\zeta,\cdot)
\end{split}
\end{equation}

{\bf Case 1: $\xi=k\eta_0$.} In this case, only the term corresponding to $\zeta=0$ survives and the differential equation reduces to
$$\partial_t \hat{u}(k\eta_0,\cdot)=-2\pi i (\hat{u}(0,\cdot)\cdot k\eta_0)\hat{u}(k\eta_0,\cdot)=-2\pi ik \hat{u}(k\eta_0,\cdot)$$
which is satisfied, because $\hat{u}(k\eta_0,t)=e^{-2\pi i kt}\hat{u}(k\eta_0,0)$.\\

{\bf Case 2: $\xi=k\eta_0+m\xi_1$, $m\neq 0$.} In this case, the differential equation reads
\begin{equation}
\begin{split}
 v_{k,m}'
=&-2\pi i\left\{e^{-2\pi i t}(\xi_1\cdot\xi) v_{k-1,m}+e^{2\pi i t}(\xi_1\cdot\xi) v_{k+1,m}+(\xi_0\cdot\xi)v_{k,m}\right\}\\
=&-2\pi i\left\{me^{-2\pi i t} v_{k-1,m}+me^{2\pi i t} v_{k+1,m}+kv_{k,m}\right\}
\end{split}
\end{equation}
which is satisfied in view of the recurrence relations \eqref{R1},\eqref{R2}.

\subsection{Proof of theorem \ref{3d-blowup}}
Given $f_1,f_2:\Z^3\to\R_+$ we have to construct functions $h$ and $g$ as above. 

With the choice of initial data above we have for $|\xi|>1$:
$$|\hat{u}_0(\xi)|=
\begin{cases}
|h(k)|, \quad \text{if} \quad \xi=k\eta_0, \;k\in\Z,\\
\\
0, \quad \text{otherwise.}
\end{cases}$$

Hence, we simply choose any $h$ such that $|h(k)|\leq f_1(k\eta_0)$.

Moreover, $|\hat{u}(m\xi_1,T+1)|=e^{-1}|g(m)|$, and we may choose $g$ such that $|g(m)|=ef_2(m\xi_1)$.

\subsection{Blow-up of Sobolev norms}\label{sec-blowup}
Now let us consider $\hat{u}(T+1,\xi)$. For $\xi=m\xi_1$ we have
$$\|\hat{u}(T+1,\xi)\|=|g(m)|/e.$$

Thus, for any $s\in\R$ we have the lower bound $$\|\hat{u}(\cdot,T+1)\|_s=\sum_{\xi\in\Z^3}(1+\|\xi\|^2)^s\|\hat{u}(\xi,T+1)\|^2\geq \frac{1}{e^2}\sum_{m\in\Z}(1+m^2)^s|g(m)|^2,$$
while, $$\|\hat{u}(\cdot,0)\|_s=\sum_{\xi\in\Z^3}(1+\|\xi\|^2)^s\|\hat{u}(\xi,0)\|^2=3+2\sum_{k\geq 2}(1+k^2)^s|h(k)|^2.$$

Hence, we have found a solution whose regularity at time $t=0$ is determined by the function $h$ and which evolves by time $t=T+1$ to a type of behaviour which is at least as bad as the decay of the function $g$ which controls the decay of the Fourier modes along the $\bra \xi_1\ket$-direction. 

\subsection{Regaining regularity}\label{sec-gain}
We note that we may take $f\in C^\infty_c(\R_+)$ such that $\supp f\subset[T_1,T_2]$, where $0<T_1<T_2$. If we repeat the above construction using this choice of $f$, then we have a solution which loses regularity on the time-interval $[T_1,T_2]$ and then regains it for any $t>T_2$. In particular, the multifractal scaling which we observed in the previous section will appear on the time-interval $[T_1,T_2]$ and disappear for $t>T_2$. If we take $f$ such that $\supp f\subset\bigcup_{j=1}^{+\infty}[T_{j,1},T_{j,2}]$, where the time-intervals $[T_{j,1},T_{j,2}]\subset\R_+$ are disjoint, then the solution loses and regains regularity infinitely many times. This coincides with the appearance and disappearance of multifractal scaling.

\subsection{Proof of Corollary \ref{mult}}\label{sec-mult-proof}
We may now rewrite the moment sum as $$\sum_{|m|\leq N}\mu_N(m,t)^q=\frac{M_{q,N}(t)}{M_{1,N}(t)^q}$$ where we introduce
$$M_{q,N}(t):=\sum_{m\in\Z}\|\hat{U}_N(m,t)\|^{2q}.$$

We have $\hat{u}(m\xi_1,t)=f(t)g(m)v$, for $m\neq 0$, and $\hat{u}(0,t)=\xi_0$. 


We write $M_{q,N}(t)$ in terms of the Fourier coefficients as
\begin{equation}
\begin{split}
M_{q,N}(t)=&\sum_{m\in\Z}\|\hat{U}_N(m,t)\|^{2q}\\
=&\sum_{|m|\leq N}\|\hat{u}(m\xi_1,t)\|^{2q}\\
=&1+|f(t)|^{2q}\sum_{0<|m|\leq N}|g(m)|^{2q}.
\end{split}
\end{equation}

Let us assume $t>T$ to ensure $f(t)\neq0$. The Renyi entropy can then be written as 
$$H_{q,N}(t)=\frac{\log(1+\sum_{0<|m|\leq N}|g(m)|^{2q})-q\log(1+\sum_{0<|m|\leq N}|g(m)|^{2})}{1-q},$$
where $q>1$. 

Let us take $g(k)=k^{-\alpha}$ for $\alpha\in(0,\tfrac{1}{2})$ and let us further assume $q>\frac{1}{2\alpha}$. 
In this case, we have, as $N\to+\infty$,
$$\sum_{0<|k|\leq N}|g(k)|^{2}\sim C_\alpha N^{1-2\alpha},$$
and $$\sum_{0<|k|\leq N}|g(k)|^{2q}=O(1).$$
Therefore, $$H_{q,N}(t)\sim (1-2\alpha)\frac{q}{q-1}\log N, \quad\text{as}\; N\to+\infty.$$

The fractal exponent equals, for any $q>\tfrac{1}{2\alpha}$,
$$
D_q(t)=\lim_{N\to+\infty}\frac{H_{q,N}(t)}{\log N}=(1-2\alpha)\frac{q}{q-1}. 
$$

\begin{remark}
The construction presented here may be applied to more general non-linear PDE. In particular, we can deal with a generalization of the system \eqref{Fourier-system}, where we add a term $\sigma(\xi)\hat{u}_j(\xi,\cdot)$, and $\sigma:\Z^3\to\R$. A detailed account will be given in \cite{U}.
\end{remark}

\section{Loss of regularity for complex data on the $2$-torus}\label{sec-complex}
Fix a vector $v\in\Z^2$. Denote $S=\bra v\ket^\perp$. Let us assume that there is a solution with the property that for all $t\geq 0$ we have $\supp\hat{u}(\xi,t)\subset S$. Moreover, let us suppose that for all $t\geq 0$ we have $\hat{u}(0,t)=i\xi_0$ for some $\xi_0\in S$. And for all $\xi\in S\setminus\{0\}: \hat{u}(\xi,t)\in \bra v\ket$.

Let us look at the r.h.s. of the first equation in \eqref{Fourier-system}. 

Suppose $\xi\notin S$: Then for each term in the sum (a solution of $\zeta+\eta=\xi$) either $\zeta\notin S$ or $\eta\notin S$, which means that that either $\hat{u}(\zeta,t)=0$ or $\hat{u}(\eta,\cdot)=0$. So all terms vanish. 
We get the differential equation $\partial_t \hat{u}_j(\xi,\cdot)=0$, and thus $\hat{u}(\xi,t)=\hat{u}(\xi,0)=0$, as $\xi\notin S$.

Suppose $\xi\in S$: Suppose that $(\zeta,\eta)\in S^2$ is a solution of $\zeta+\eta=\xi$ (note that all other terms vanish). If $\zeta\neq 0$, then $\hat{u}(\zeta,t)\in \bra v\ket$ and therefore $\hat{u}(\zeta,t)\cdot \eta=0$. Hence the only nonzero term in the sum is the one corresponding to the solution $\zeta=0$, $\eta=\xi$. In this case we obtain the differential equation
\beq
\partial_t \hat{u}_j(\xi,\cdot)=-2\pi i(\hat{u}(\xi,\cdot)\cdot \xi)\hat{u}(\xi,\cdot)=2\pi(\xi_0\cdot\xi)\hat{u}(\xi,\cdot)
\eeq
We solve to get $\hat{u}(\xi,t)=e^{2\pi(\xi_0\cdot\xi)t}\hat{u}(\xi,0)$. In particular $\hat{u}(0,t)=\hat{u}(0,0)=i\xi_0$.

So to construct a solution we take a function $f:S\setminus\{0\}\mapsto \C$. 

Then take 
\beq
\hat{u}(\xi,0)=
\begin{cases}
i\xi_0, \quad\text{if}\quad \xi=0,\\
\\
f(\xi)v, \quad\text{if}\quad \xi\in S\setminus\{0\},\\
\\
0, \quad\text{if}\quad \xi\notin S.
\end{cases}
\eeq
As a solution we get
\beq    
\hat{u}(\xi,t)=e^{2\pi(\xi\cdot\xi_0)t}\hat{u}(\xi,0).
\eeq
We verify that this solution satisfies all of our assumptions.

We compute
$$\int_{\T^2}|u(x,t)|^2 dx = |\xi_0|^2+\|v\|^2\sum_{\xi\in S\setminus\{0\}} |f(\xi)|^2e^{4\pi(\xi\cdot\xi_0)t}.$$

Let us choose an initial condition corresponding to the function $f(\xi)=g(\xi)e^{-\gamma|\xi|}$ on the hyperplane $S$.

We have for $t\in[0,T]$, where $T=\gamma/(2\pi |\xi_0|)$,
\beq
\begin{split}
\|u(\cdot,t)\|_s^2&=\sum_{\xi\in\Z^2}(1+|\xi|^2)^s|\hat{u}(\xi,t)|^2\\
 &\leq |\xi_0|^2+\|v\|^2\sum_{\xi\in S\setminus\{0\}} (1+|\xi|^2)^s|g(\xi)|^2e^{(4\pi|\xi_0|t-2\gamma)|\xi|}
\end{split}
\eeq

Moreover, we have the lower bound
\beq
\begin{split}
\|u(\cdot,t)\|_s^2=&\sum_{\xi\in\Z^2}(1+|\xi|^2)^s|\hat{u}(\xi,t)|^2 \\ 
\geq& \|v\|^2\sum_{n\geq 0} (1+|n\xi_0|^2)^s|g(n\xi_0)|^2 e^{(4\pi n|\xi_0|^2 t-2\gamma n|\xi_0|)}
\end{split}
\end{equation} 

Let us suppose that $|g(\xi)|\geq |\xi|^{-\alpha}$ for $\xi=k\xi_0$, $k\geq1$, for certain $\alpha>1/2$.

So we see that 
\beq
\|u(\cdot,t)\|_s
\begin{cases}
<+\infty, \quad \text{for}\; t<T \quad \text{and for all}\; s\in\R\\
=+\infty, \quad \text{for}\; t=T \quad \text{and for all}\; s\geq\alpha-\tfrac{1}{2}\\
=+\infty, \quad \text{for}\; t>T \quad \text{and for all}\; s\in\R
\end{cases}
\eeq


\begin{thebibliography}{99}
\bibitem{dLSz09}
C. de Lellis, L. Sz\'ekelyhidi Jr., {\em The Euler equations as a differential inclusion}, Ann. Math 170 (3) (2009), pp. 1417--36.
\bibitem{dLSz13}
C. de Lellis, L. Sz\'ekelyhidi Jr., {\em Dissipative continuous Euler flows}, Invent. Math. 193, Issue 2 (2013), Page 377-407.
\bibitem{FFT75}
C. Foias, U. Frisch, R. Temam, {\em Existence des solutions $C^\infty$ des \'equations d'Euler}, C. R. Acad. Sc., Paris, 280 A (1975), pp. 505--508.
\bibitem{FT89}
C. Foias, R. Temam, {\em Gevrey Class Regularity for the Solutions of the Navier-Stokes Equation}, J. Funct. Anal. 87 (1989), 359--369.
\bibitem{SM99}
J. C. Mattingly, Ya. G. Sinai, {\em An Elementary Proof of the Uniqueness and Existence Theorem for the Navier-Stokes Equations}, Comm. Contemp. Math. 1 (1999), 497--516.
\bibitem{CS17} D. Li, Ya. G. Sinai, {\em Blow-ups of complex solutions of the 3D Navier-Stokes system and renormalization group method}, J. Eur. Math. Soc. 10 (2008), no. 2, 267--313.
\bibitem{Sch93}
V. Scheffer, {\em An inviscid flow with compact support in space-time}. J. Geom. Analysis 3 (4) (1993), pp. 343--401.
\bibitem{Sn96}
A. Snirelman, {\em On the non-uniqueness of weak solution of the Euler equation}, Journ\'ees E. D. P. (1996), pp. 1--10.
\bibitem{KU1}
J. P. Keating, H. Uebersch\"ar, {\em Multifractal eigenfunctions for a singular quantum billiard}, Comm. Math. Phys.  Vol. 389 (2022), 543--569
\bibitem{KU2}
J. P. Keating, H. Uebersch\"ar, {\em Multifractal eigenfunctions for quantum star graphs}, arXiv:2202.13634.
\bibitem{U}
H. Uebersch\"ar, {\em Multifractal scaling and periodic solutions to non-linear PDE}, in prep.
\end{thebibliography}
\end{document}